\documentstyle{amsppt}
\magnification\magstep1
\NoRunningHeads
\hsize6.9truein
\pageheight{23 truecm}
\baselineskip=15pt
\catcode`@=11
\def\logo@{\relax}
\catcode`@=\active


\def\t#1{\tilde{#1}}

\def\ideal#1.{I_{#1}}
\def\ring#1.{\Cal O_{#1}}
\def\proj#1.{\Bbb P(#1)}
\def\pr #1.{\Bbb P^{#1}}
\def\af #1.{\Bbb A^{#1}}
\def\Hz #1.{\Bbb F_{#1}}
\def\Hbz #1.{\overline{\Bbb F}_{#1}}
\def\fb#1.{\underset #1 \to \times}
\def\res#1.{\underset {\ \ring #1.} \to \otimes}
\def\au#1.{\operatorname {Aut}\,(#1)}
\def\deg#1.{\operatorname {deg } (#1)}
\def\pic#1.{\operatorname {Pic}\,(#1)}
\def\pico#1.{\operatorname{Pic}^0(#1)}
\def\ner#1.{NS (#1)}
\define\rdown#1.{\llcorner#1\lrcorner}\
\define\rup#1.{\ulcorner#1\urcorner}


\def\list#1.#2.{{#1}_1,{#1}_2,\dots,{#1}_{#2}}
\def\loc#1.#2.{\Cal O_{#1,#2}}
\def\fderiv#1.#2.{\frac {\partial #1}{\partial #2}}
\def\map#1.#2.{#1 \longrightarrow #2}
\def\rmap#1.#2.{#1 \dasharrow #2}
\def\emb#1.#2.{#1 \hookrightarrow #2}
\def\non#1.#2.{\text {Spec }#1[\epsilon]/(\epsilon)^{#2}}
\def\H#1.#2.{\text {Hilb}^{#1}(#2)}
\define\sym#1.#2.{\operatorname {Sym}^{#1}(#2)}
\def\Hb#1.#2.{\text {Hilb}_{#1}(#2)}
\def\Hm#1.#2.{\Hom_{#1}(#2)}
\def\prd#1.#2.{{#1}_1\cdot {#1}_2\cdots {#1}_{#2}}


\def\alist#1.#2.#3.{#1_1\,#2\,#1_2#2\,\dots\,#1_{#3}}
\def\lmap#1.#2.#3.{#1 \overset #2\to \longrightarrow #3}
\def\ses#1.#2.#3.{0\longrightarrow #1 \longrightarrow #2 \longrightarrow #3 
\longrightarrow 0}
\def\les#1.#2.#3.{0\longrightarrow #1 \longrightarrow #2 \longrightarrow #3}

\def\Hi#1.#2.#3.{\text {Hilb}^{#1}_{#2}(#3)}


\define\Hom{\operatorname{Hom}}

\def\deg{\operatorname{deg}}

\def\sg{\operatorname{Sing}}


\def\mm{\overline{M}}
\def\lm{\widetilde{M}}
\def\mn{\overline{M}_{0,n}}
\def\ln{\widetilde{M}_{0,n}}
\def\mg{\overline {\Cal M}_g}


\def\mapdown#1{\big\downarrow\rlap{$\vcenter
{\hbox{$\scriptstyle#1$}}$}}

\def\mapse#1{
{\vcenter{\hbox{$\mathop{\smash{\raise1pt\hbox{$\diagdown$}\!\lower7pt
\hbox{$\searrow$}}\vphantom{p}}\limits_{#1}\vphantom{\mapdown{}}$}}}}

\topmatter
\title Contractible extremal rays on $\overline{M}_{0,n}$. \endtitle
\author Se\'an Keel
 and James M\raise 1.6pt \hbox{\text {\smc c}}Kernan
 \endauthor
\address Department of Mathematics, 
University of Texas, at Austin, Austin TX 78712 
Department of Mathematics,
University of California at Santa Barbara
Santa Barbara, CA 93101
\endaddress
\endtopmatter
\def\vt{\Cal T}
\def\uu{\Cal U}
\heading \S 1 Introduction and statement of results \endheading

One of the richest objects of study in higher dimensional algebraic geometry 
is the Mori-Kleiman (closed) cone of curves, $\overline{NE_1}(M)$, defined as 
the closed convex cone in $H_2(M,{\Bbb R})$ generated by classes of 
irreducible curves on $M$. A lot of geometric information about $M$ is encoded
in the cone of curves. For example the possibilities for maps with connected 
fibres are determined by the cone's faces. Not surprisingly, 
$\overline{NE_1}(M)$ is difficult to compute. Even when $M$ is well
understood, it can be difficult to 
find generators for the cone, 
that is, to find all of the ``edges'', or to use the technical
term, ``extremal rays'' (``edge'' is potentially misleading as portions of 
the cone may be curved). Indeed, it is not even generally obvious whether or
not a given curve spans an extremal ray. One expects better luck
understanding rays 
on 
which $-c_1(M) = K_M$ (or more 
generally log terminal $K_M+\Delta$) are negative, as these 
are described by the 
powerful cone and contraction theorems of Mori-Kawamata-Shokurov: each is 
generated by a smooth rational curve, and can be ``contracted'', i.e. there 
is a map (with domain $M$) whose fibral curves are precisely the curves
whose homology class lies on the extremal ray. Thus as a first step
in computing the cone, it is natural to consider such ``negative''
extremal
rays, or more generally, rays which can be contracted.

Here we consider $\mn$, the moduli space of stable $n$-pointed rational 
curves, as well as $\ln$ the quotient of $\mn$ by the natural symmetric group 
action, which is (an irreducible component of) the moduli space of log pairs 
(see \cite{Alexeev94}). 

The locus of points in $\mn$ corresponding to a curve with at least $k+1$ 
components has pure codimension $k$; we call its irreducible components the 
{\bf vital codimension $k$-cycles}. Vital divisors, curves, $k$-cycles etc. 
are analogously defined. By a vital cycle in $\ln$ we mean the image of a 
vital cycle in $\mn$. It is relatively easy to check that the vital cycles 
generate the Chow group. It is natural to wonder if much more is true:

\proclaim{1.1 Question} Is every effective cycle linearly
equivalent to an effective sum of vital cycles?
\endproclaim

This was first posed to us by William Fulton. In the interest of drama, we 
will refer to (1.1) as Fulton's conjecture. Here we consider only the cases 
of curves and divisors. As homological and linear equivalence are the same 
on $\mn$, the conjecture in these cases is equivalent to the statement that 
vital cycles generate all extremal rays of $\overline{NE}_1$ and 
$\overline{NE}^1$, the cones of curves and divisors. We prove this for 
$\overline{NE}^1(\ln)$ and for contractible extremal rays of 
$\overline{NE}_1(\mn)$. 

Let $D \subset \mn $ be the boundary, i.e. the sum of the vital divisors. 
Let $D = \sum B_i$ be its decomposition into $S_n$ orbits (there are $[n/2]$ 
such orbits). For a subvariety $Z \subset \mn $, let $\t{Z}$ be its image 
(with reduced structure) in $\ln$. 

Here are precise statements of our results:

\proclaim{1.2 Theorem} Let $R$ be an extremal ray of the cone of curves 
$\overline{NE}_1(\mn)$. Then $R$ is spanned by a vital curve under any of 
the following conditions
\roster
\item There is a morphism  $f:\map \mn.Y.$, contracting $R$, 
with $\rho(Y) = \rho(\mn) -1$, and such that the exceptional locus of $f$ 
is not a curve. 
\item $(K_{\mn}+ G)\cdot R < 0$, where $G$ is an effective boundary whose
support is contained in $D$.
\item $n \leq 7$. 
\endroster
\endproclaim

Of course (1.2.3) says Fulton's conjecture holds for curves, 
provided $n\leq 7$. We were able to prove much stronger results for 
$\ln$ (especially (1.3.1-2)):
\proclaim{1.3 Theorem} 
\roster
\item The cone of effective divisors $NE^1(\ln)$ 
is simplicial, generated by the $\t{B}_i$. 
\item An effective divisor on $\ln$ fails to be big iff its support is a 
proper subset of $\t{D}$. Any non-trivial nef divisor is 
big. 
\item The cone of curves of $\overline{NE}_1(\ln)$ is generated by curves in 
$\t{D}$.
\endroster

Now suppose $n \leq 11$.
\roster
\item[4] $NE_1(\ln)$ is a finite rational polyhedron, with edges spanned by 
images of vital curves, 
\item Every proper face is contractible by a log Mori fibre space. In 
particular every nef divisor is eventually free. 
\item The divisor $\sum_{i=2}^{[n/2]} r_i \t{B}_i$
is nef (resp. ample) iff 
$$
r_{a+b} + r_{a+c} + r_{b + c} - r_a -r_b - r_c - r_d 
$$
is non-negative (resp. strictly positive), for all positive integers
$a$, $b$, $c$ and $d$, with $n=a+b+c+d$ (where we define $r_1=0$ and 
$r_i=r_{n-i}$ for $i > [n/2]$).
\endroster
\endproclaim

(By a simplicial cone we mean a cone over a simplex, i.e. a polyhedral
cone whose edges are linearly independent) 

The spaces $\mn$ and $\ln$ are interesting from a number of viewpoints. 
They are closely related to the moduli space of curves, $\mg$.
A finite quotient of $\mn$ occurs as a locus of degenerate curves in the 
boundary of $\mg$, while $\mn$ is the base of the complete Hurwitz scheme 
(see \cite{HarrisMumford82}) 
which can be used, for example, to prove that $\mg$
is irreducible. By \cite{Kapranov93a}, $\mn$ parameterises degenerations of 
rational normal curves. Generalisations of $\mn$ are important for Quantum 
Cohomology calculations, see \cite{KonsevichManin94}. 
$\mn$ is useful for studying 
fibrations with general fibre $\pr 1.$, as in particular it can sometimes 
be used in lieu of a minimal model program. Kawamata exploits this in 
\cite{Kawamata77} to prove additivity of log Kodaira dimension for one 
dimensional fibres, and in \cite{Kawamata95} to prove a codimension two 
subadjunction formula. 

We note that there is an explicit construction of $\mn$, as a 
blow up of $\pr {n-3}.$ along a sequence of simple centres (see (3.1)).
In particular $\mm_{0,5}$ is a del Pezzo of degree five, $\mm_{0,6}$ is log 
Fano, and $\mm_{0,7}$ is nearly log Fano, in the sense that
$-K_{\mm_{0,7}}$ is effective.  We do not know of such an explicit 
contruction of $\ln$, and we have in general a much weaker grasp on its
geometry (though a much stronger grasp on its cones). Note by (1.3.3),
$\ln$ admits no nontrivial fibrations. See also (3.7). 

Despite the fact that the blowup construction gives an easy computation 
of some invariants of $\mn$, one cannot expect the same for the cone of 
curves. For example the blow up of $\pr 2.$ in
eight points has a finite polyhedral cone of curves, but one can
choose a ninth point in such a way that the blow up has
a cone with infinitely many edges. We do not use the blow up
description in any significant way in our proof of (1.2-3).

As we note in (3.5) $-K_{\ln}$ and $-K_{\mn}$ are not effective for 
$n \geq 8$. In view of this, the cases of (1.3) for $8 \leq n \leq 11$ 
are interesting in that they give examples of non log Fano varieties, 
for which every face of the cone of curves is none the less contractible. 
From this perspective, (1.1), if true, would really be rather 
surprising. Each vital curve (indeed every vital cycle) is smooth and
rational. The cone of curves of a log Fano is generated by
rational curves, but one does not expect this in general, even for a 
rational variety. For example, let $S$ be the blow up of $\pr 2.$ in a large 
number of general points. As observed by Koll\'ar, and independently by 
Caporaso and Harris, $K_S$ is strictly negative on rational curves, but of 
course $K_S^2 < 0$, so $K_S$ must be positive on some curves (but see the 
remark after (2.4)). 

If (1.1) holds for curves, then one can describe the
ample cones of $\mn$ or $\ln$ by a series of inequalities analogous
to those in (1.3), using (4.3). One can then describe, at
least in theory, the cone of curves, since it is dual to the ample cone.
As an example of the complexity of these cones, 
$NE_1(\mm _{0,7})$ is a polyhedral cone of dimension $42$ with $350$ edges
(see (4.3) and (4.6)). 

Fulton's conjecture implies every vital curve spans an extremal ray
and each is $K_{\mn}+ G$ negative for some $G$ as in (1.2.2) (see (4.6)).
So by the contraction theorem \cite{KMM87} each vital curve
is contracted by a map of relative Picard number one. For $n \geq 9$ every
vital curve deforms. So if (1.1) holds, then (1.2) contains all 
the possibilities for extremal rays, and 
(1.2.1) has all the possibilities for $n \geq 9$.

Here is a brief outline of our proofs of (1.2-3). It turns out each component
of $D$ is a product $\mm _{0,i} \times \mm _{0,j}$, for $i$, $j < n$ 
(see \S 3). For (1.2) we proceed by induction, the main work is to show the 
extremal ray $R$ is in the subcone generated by curves in $D$. For this
our main tool is (2.2). (1.3) follows from (1.2) and
some simple intersection calculations: one set to show that $NE^1(\ln)$ is
simplicial, and a second to show that for $n \leq 11$ every face of
the cone contracts to a log Mori fibre space. 

\S 2 contains some results about the cone spanned by curves lying in a 
divisor. Most of the results of \S 2 are of independent interest
(in particular (2.3-2.6)), and 
hopefully have broader applications.
For this reason we work in greater generality. 
\S 3 contains the necessary 
ingredients to apply some of the results of \S 2 to $D \subset \mn $. 
Intersection products of various vital cycles are easy to compute, and the 
pairing between divisors and curves is described in \S 4. \S 5  
finishes the proof of (1.3). 

We would like to say a few words about other seemingly natural approaches to
(1.1). For curves, it is enough (in fact equivalent) to show that
if a divisor intersects all vital curves non-negatively, then it
is nef. By induction it is sufficient to show that such a divisor 
is linearly equivalent to an effective sum of vital divisors.
As the vital divisors generate the Picard group, the intersection
conditions give a finite collection of simple inequalities on
the coefficients. Unfortunately the 
combinatorics are intimidating, and we were not able to make any progress
in this direction, even for $n=6$.

Throughout we will use the main results of the minimal model
program, the contraction theorem, the cone theorem etc., as well as
the established notation as set out in \cite{KMM87}.
We also use elementary properties  and notions
of cones from Chapter II.4 of \cite{Koll\'ar96}. In particular,
by an extremal ray $R$ of a closed convex cone $W$ we mean
a one dimensional subcone with the property that if
$x + y \in R$ for $x,y \in W$ then $x,y \in R$. We note
that every closed convex cone is the convex hull of its
extremal rays. All spaces are assume to be of finite type over $\Bbb C$. 
Unless otherwise stated, by a divisor we mean an $\Bbb R$-divisor. In
\cite{Shokurov96} the main results of the MMP are extended to
$\Bbb R$-divisors. However we only need one such result (see (2.2)). 

\heading \S 2 The cone spanned by curves inside a divisor \endheading

 We first introduce some notation and definitions. Let $D$ be a reduced 
Weil divisor inside the projective $\Bbb Q$-factorial klt variety $M$ of 
dimension $n$. Let $W$ be the closed subcone of $\overline{NE}_1(M)$ 
generated by curves lying in $D$. 

We are interested in extremal rays that lie outside of $W$ and moreover under 
what conditions $W=\overline{NE}_1(M)$. 

\definition{2.1 Definitions}  We say that $D$ has 
{\bf anti-nef normal bundle} if for every curve 
$C\subset D$, $C\cdot D\leq 0$. 

We will say an extremal ray $R$ of $\overline{NE}_1(M)$ is {\bf log extremal}
if there exists a klt divisor $K_M+\Delta$ such that $(K_M+\Delta)\cdot R<0$.
\enddefinition

Log extremal rays are very special: By the cone and
contraction Theorems they are spanned by
rational curves $C$, and there is a morphism $f:\map M.Y.$ contracting $C$ 
such that $f_*(\ring M.)=\ring Y.$ and $\rho(Y)=\rho (M)-1$. 

 The following recent result of Shokurov will prove useful:

\proclaim{2.2 Lemma(Shokurov)} Let $X$ be a projective variety, and let 
$L \in N^1(X)$ be a nef class (not necessarily rational) 
with $L^{dim(X)} > 0$. Then $L$ is in the interior of $NE^1(X)$. 
\endproclaim
\demo{Proof} This is implied by the proof of (6.17) of \cite{Shokurov96}.
\qed\enddemo

(2.2) has some interesting corollaries:

\proclaim{2.3 Corollary}  If the components of $D$ span $\overline{NE}^1(M)$
then $W=\overline{NE}_1(M)$.
\endproclaim
\demo{Proof} Let $D = \sum D_i$ be the decomposition of $D$ into irreducible 
components.

Let $A$ be an ample divisor with support in $D$, and
let $R \subset \overline{NE}_1(M)$ be an extremal ray. Assume $R \not \in W$.
Let $L$ be a nef class supporting $R$. $L|D$ is ample. Since 
$L$ is an effective sum of $D_i$, $L^{\dim M} > 0$ thus
by (2.2), $R$ cannot be numerically effective. Since
the $D_i$ generate $\overline{NE}^1(M)$, $R \cdot D_i <0$
for some $i$. But this implies $R \in W$, a contradiction. \qed\enddemo

\proclaim{2.4 Proposition} Let $G$ be an  effective $\Bbb Q$-divisor, with 
non-empty support $D$. 

 Let $R$ be an extremal ray of $\overline{NE}_1(M)$, which does not lie in 
$W$. If $(K_M+G)\cdot R\leq 0$ then $R$ is log extremal and 
$K_M\cdot R\leq 0$. 

 In particular, if $-(K_M+G)$ is nef then $\overline{NE}_1(M)$ is spanned 
by $W$ and log extremal rays $R$, such that $K_M\cdot R\leq 0$. 
\endproclaim
\demo{Proof} Let $R$ be an extremal ray of $\overline{NE}_1(M)$, not lying 
in $W$. In particular $R \cdot D_i\geq 0$ and so $K_M\cdot R\leq 0$.

 On the other hand we are done if $K_M\cdot R<0$. Thus we may assume 
$K_M\cdot R=0$. Let $L \in N^1(M)$ be a nef class supporting $R$. Then $L$ 
is strictly positive on $W\setminus 0$ and so by compactness of a slice of 
$W$, $L+\epsilon D$ is nef and supports $R$ for $0< \epsilon \ll 1$. As $L$ 
is ample on $D$ $L^{n-1} \cdot D > 0$. In particular we can replace $L$ by 
$L+\epsilon D$ and assume $L^n> 0$. Then by (2.2) $R \cdot V < 0$ for some 
effective Weil divisor $V$. But $(K_M+\epsilon V)\cdot R <0$  and
$(K_M+\epsilon V)$ is klt for $0<\epsilon \ll 1$. \qed\enddemo

\remark{Remark} (2.4) is interesting even in the case of a surface. For 
example pick a cubic in $\pr 2.$ and blow up as many points as
you like along the cubic. Let $M$ be the resulting surface and $D$ the 
strict transform of the cubic. (2.4) then says that $D$ union all the 
$-2$-curves and $-1$-curves generate the cone of curves of $M$. 
\endremark

\proclaim{2.5 Proposition} Let $f:\map M.Y.$ be a proper surjection from a 
smooth projective variety $M$ to a normal variety $Y$ with 
$f_*(\ring M.)=\ring Y.$, and $\rho(Y) = \rho(M)-1$. Suppose $D$ has ample
support and each irreducible component of $D$ has anti-nef normal bundle. 

If $f|D$ is finite then $f$ is birational, and its exceptional locus is 
a curve. 
\endproclaim
\demo{Proof} Suppose on the contrary that there is an irreducible 
surface $E$ whose image has dimension at most one.

 Let $D= \sum _iD_i$ be the decomposition of $D$ into irreducible components.
Note the assumptions on Picard number imply that any class in $N^1(M)$ which 
is zero on some fibral curve, is pulled back from $N^1(Y)$. 

Since $D$ has ample support $I=D\cap E$ is non-empty. As $f|D$ is finite, 
$I$ and each $D_i \cap E = D_i \cap I$, is an effective $\Bbb Q$-Cartier 
divisor of $E$, and in particular, is purely one dimensional. Thus
if $I$ meets $D_i$, it has an irreducible component contained
in $D_i$. Since $D$ has ample support, and $f|D$ is finite,
$E$ contracts to an irreducible curve $C \subset f(D)=f(I)$ and $f|I$ is 
finite.\newline 
\hskip .6truecm {\it Claim.} \hskip .1truecm We can find two irreducible 
components $B_1$, $B_2$ of $I$ and (after renaming) two divisors $D_1$, $D_2$ 
with $B_i \subset D_i$ such that $B_i \cdot D_j \geq 0$ (for $i \neq j$) 
and at least one inequality is strict:

 Choose an irreducible component $B_1$ of $I$ contained in a maximal number 
of $D_i$. Suppose (after reordering) $\list D.k.$ are the components
of $D$ containing $B_1$. Since the $D_i$ have anti-nef
normal bundles, and $D$ has ample support, for some $j > k$ we have
$D_j \cdot B_1 > 0$. Let $B_2$ be an irreducible component
of $D_j \cap I$.  By the choice of $B_1$
we can assume (after reordering) that $B_2 \not \subset D_1$.
Now set $D_2= D_j$. 

This establishes the claim. 

Since $D_1$, $D_2$ each meet a fibre, we can choose $\lambda > 0$ such that
$D_1 - \lambda D_2$ is pulled back from $Y$. Let $J = (D_1 - \lambda D_2)|E$. 
Then $J \cdot B_1 \leq 0$ and $J \cdot B_2 \geq 0$, and one inequality is 
strict. Since $J$ is pulled back from $C$, and the $B_i$ are multi-sections, 
this is a contradiction. \qed\enddemo

\remark{Remarks} 
\roster
\item The assumption on the relative Picard number in (2.5) is necessary;
it cannot be replaced by the weaker assumption that $f$ is the contraction
of an extremal ray. For example consider $M=E\times E$ for an elliptic curve 
$E$, $D = F_1 + F_2$ the sum of the two fibres and $f:\map M.E.$ the addition 
map. 
\item The assumption on Picard number holds 
when $f$ is the contraction of a log extremal ray.
\item One can not rule out the final possibility. For example:
Let $M$ be a del Pezzo surface whose cone of curves
is not a simplex (e.g. blow up $\pr 2.$ at three non-collinear points). 
Let $D$ be a sum of $\rho(M)$ $-1$-curves with ample support
(any effective class is a sum of at most $\rho(M)$ extremal rays,
and all the extremal rays are $-1$-curves). Let $f$ blow down
some other $-1$-curve.
\endroster
\endremark

\proclaim{2.6 Lemma} Suppose $M$ is smooth of dimension at least
three, every component of $D$ has 
anti-nef normal bundle, and $D$ has ample support. Let $G$ be a 
nonempty effective 
$\Bbb Q$-divisor whose support lies in $D$. 
\roster 
\item Let $R$ be an extremal ray of $\overline{NE}_1(M)$. If either 
$(K_M+G)\cdot R<0$, or the dimension
of $M$ is at least four and $(K_M+G)\cdot R\leq 0$ then $R\in W$.
\item If $-(K_M+G)$ is nef, and either the support of $G$ is exactly $D$ or
the dimension of $M$ is at least four, then $W=\overline{NE}_1(M)$. 
\endroster
\endproclaim
\demo{Proof} Let $R$ be an extremal ray of $\overline{NE}_1(M)$, and 
suppose $R\notin W$ but $(K_M+G)\cdot R\leq 0$. Then by (2.4) we know
that $R$ is spanned by a contractible rational curve $C$. (1) and (2) 
now follow easily from (2.5) and the observation that if 
$K_M\cdot C<0$ and $M$ is a threefold (resp. $K_M\cdot C\leq 0$ and $M$ 
has dimension at least four) then $C$ deforms inside $M$ 
(see II.1.13 of \cite{Koll\'ar96}). 
\qed \enddemo
We will use the following technical result in the next section. 

\proclaim{2.7 Lemma} Let $N\subset M$ be a normal divisor 
and suppose that $\map {N^1(M)}. N^1(N).$ is surjective.
Let $f:\map M.Y.$ be a map to a normal projective variety with
$f_*(\ring M.) = \ring Y.$, and $\rho(Y) = \rho(M) -1$. Let $g:\map N.Z.$
be the Stein factorisation of $f|_N$. If $f|N$ is not 
finite, then $\rho(Z) = \rho(N) -1$.
\endproclaim
\demo{Proof} $f$ contracts an extremal ray $R$. Suppose $f|N$ is
not finite. Then $R \in N_1(N)$. If $L \in N^1$ and 
$L \cdot R =0$, then $L|N$ is pulled back from $Z$. Since every class in 
$N^1(D)$ extends to $M$, the result follows. \qed \enddemo

\heading \S 3  Geometry of $\mn$ and $\ln$.  \endheading

We will use (a slight modification of) the notation of, as well as several 
simple facts from pg. 551--554 of \cite{Keel92}. For the readers convenience 
we will recall the most important ideas:

A vital divisor is determined by a partition of $\{1,2,\dots,n\}$
into disjoint subsets $T$, $T^c$, each containing at least two
elements. The generic point of the corresponding vital divisor $D_{T,T^c}$ 
is a curve with two irreducible components, with the labels of $T$ on one 
component, and the labels of $T^c$ on the other. There is a canonical
isomorphism
$$
D_{T,T^c} =M_{T \cup \{b\}} \times M_{T^c \cup \{b\}}
$$
where e.g. by $M_{T \cup \{b\}}$ we mean a copy of $\mm_{0,|T| + 1}$ 
with the indices labeled by the elements of $T$, with $b$ an
extra index, corresponding to the singular point. We indicate
the two projections by $\pi_T$ and $\pi_{T^c}$.

The vital divisors have normal crossings, and each vital codimension
$k$-cycle is uniquely expressible as a complete intersection of
vital divisors. Each vital $k$-cycle has an expression as a
product of $\mm_{0,i}$ analogous to that for the vital divisors. In
particular, under the above decomposition, any vital curve of $D_{T,T^c}$ is 
a product of a vital curve on one factor, with a vital point on the 
second.

\proclaim{3.1 Proposition(Kapranov)} For each index $i \in \{1,2,\dots,n\}$ 
there is a birational map $q_i:\map \mn.{\pr n-3.}.$ with the following
properties:
\roster
\item $q_i$ is a composition of blow ups along smooth centres, constructed as 
follows. Fix $n-1$ general points, and blow up successively (from lowest
to highest dimensional) the (strict transforms) of every 
linear subspace spanned by any subset of these points.  
\item $q_i$ takes vital cycles to to linear spaces spanned by the chosen 
points.
\item If $i \in T$ then $q_i|_{D_{T,T^c}} = q_i \circ \pi_T$ for $i \in T$. 
\item If $F$ is the general fibre of the map $\map \mn.\mm_{0,n-1}.$ given by
dropping the $i^{th}$ point, then $q_i(F)$ is a rational normal
curve. 
\item $q_i$ is a composition of smooth blow downs, blowing down
iteratively the (images of) the divisors $D_{T,T^c}$ with 
$i \not \in T$, and $|T|=3,4,\dots,n-2$.
\endroster
\endproclaim
\demo{Proof} See \cite{Kapranov93b}.\qed \enddemo

\proclaim{3.2 Lemma} Let $\phi$ be an element of $\au \pr 1..$ of
finite order $p$ and let $Z$ be the closure of the locus of 
$n$-tuples of distinct elements of $\pr 1.$ whose coordinates
are permuted by $\phi$. Let $q$ be a general element of $Z$. If
$\phi$ fixes $j$ coordinates of $q$ then the dimension of $Z$
at $q$ is $(n-j)/p$. \endproclaim
\demo{Proof}
 Let $G \subset \au \pr 1. .$ be the subgroup
generated by $\phi$. Then $G$ has a non-trivial finite orbit,
from which it follows that $G$ has exactly two fixed points,
and after changing coordinates (so the fixed points are
$0$ and $\infty$) $\phi : \af 1. \rightarrow \af 1.$ is multiplication
by a $p^{\text{th}}$ root of unity. The coordinates divide into
orbits, each of which is either a fixed point, or has exactly
$p$ elements. The result follows. \qed  \enddemo

\proclaim{3.3 Lemma} $S_n$ acts freely in codimension one on 
$M_n \setminus D$ for $n \geq 7$, and faithfully for 
$n \geq 5$.  

The action of $S_4$ on $M_4$ factors
through the action on the set of partitions of $\{1,2,3,4\}$ into
disjoint subsets of two elements. Nontrivial elements of the 
kernel are of form $(ij)(kl)$ for $i,j,k,l$ distinct.
\endproclaim
\demo{Proof} The claims about the $S_4$ action are easily
checked, and left to the reader. 

So assume $n \geq 5$. We bound the dimension of the locus
of points in $M_n \setminus D$ 
which are fixed by some element of $S_n$. Equivalently,
we bound the dimension of the locus
$Z \subset M_n \setminus D$ of $n$-tuples (modulo automorphisms)
whose coordinates are permuted by some automorphism of $\pr 1.$.

Let $\uu \subset {\pr 1.}^{\times n}$ be the locus of distinct
points, and let 
$$
\vt = \bigl \{ (q,\phi,a,b) \bigl | q \in \uu, \phi \in Aut(\pr 1.) 
a \neq b \in \pr 1. \text{ s.t. }
\phi \text{ permutes } q \text{ and fixes }
(a,b) \bigl \}.
$$

Replace $\vt$ by any one of its irreducible components.

Let $\phi$ be a general point of the image of 
$pr: \vt \rightarrow Aut(\pr 1.)$, and let $q \in pr^{-1}(\phi)$.
By (3.2) $pr^{-1}(\phi)$ has dimension $(n-j)/p$ at $q$, while the fibre
of $\vt \rightarrow U$ has dimension three. Thus at the image of
$q$, 
$U$ has dimension $(n-j)/p -1$. The result follows. \qed
\enddemo

Let $B_i = \sum_{|T| = i} D_{T,T^c}$ for $ 2 \leq i \leq k=[n/2]$.
$B_i$ is the orbit under $S_n$ of any $D_{T,T^c}$ with $|T| = i$.

\proclaim{3.4 Lemma} For $n \geq 7$ the quotient map $q:\map \mn.\ln.$
is unramified in codimension one outside of $B_2$, and has ramification 
index two along $B_2$. 
\endproclaim
\demo{Proof} Suppose $\sigma \in S_n$ fixes each point of the 
irreducible divisor $G \subset \mn$. By (3.3), $G = D_{T,T^c}$ for some 
$T$ preserved by $\sigma$. Since the action of $\sigma$ on $M_{T \cup \{b\}}$ 
factors through the subgroup of $S_{|T| +1}$ which fixes $b$, it follows 
from (3.3) that $T = \{i,j\}$ and $\sigma= (i,j)$. \qed\enddemo

We will use the following formulae, essentially from
\cite{Pandhapripande95}:
\proclaim{3.5 Lemma} 
$$
K_{\mn}+\sum_{j=2}^k(2 - \frac{j(n-j)}{n-1})B_j=0=
K_{\ln}+(\frac 12+\frac 1{(n-1)})\t{B_2}+
\sum_{j=3}^k(2 - \frac{j(n-j)}{n-1})\t{B_j}
$$

In particular $-K_{\ln}$ (resp. $-K_{\mn }$) is pseudo-effective iff 
$n\leq 7$.
\endproclaim
\demo{Proof} The first formula is Proposition 1 of 
\cite{Pandhapripande95},
the second follows easily from the first and (3.4) and the last statement
then follows from (4.8). In fact we may use (4.3) to prove the first 
formula in a similar way to the way it is derived in 
\cite{Pandhapripande95}.  
 
 However it is possible to prove the first formula in an entirely elementary
way, using (3.1). Indeed the image $D'$ of $D$ is the union of $\binom {n-1}2$ 
hyperplanes, and the coefficients of $B_i$ are easily identified as the 
discrepancies of the divisor $K_{\pr n-3.}+(2/(n-1))D'$. \qed\enddemo

\proclaim{3.6 Lemma} $K_{\mn}+ D$ is ample and is linearly
equivalent to an effective divisor with the same support as
$D$.
\endproclaim
\demo{Proof} We proceed by induction on $n$. The result is
easy for $n=4$.

By (3.5), $K_{\mn}+ D$ is linearly equivalent
to an effective divisor $\Gamma$ with support $D$.

Note that $(K_{\mn}+ D)|D_T$
is the tensor product of the ``same expressions'' pulled back from
the two components in the product description of $D_T$. Thus
by induction $(K_{\mn}+ D)|D$ is ample.

It is easy to see that $D$ meets (set theoretically) every curve. 
Use induction and consider the map $f:\map \mn.\mm_{0,n-1}.$, observe
that $D$ meets every fibral curve, and note that 
$D \supset f^{-1}(D(\mm_{0,n-1}))$.

Thus $(K_{\mn}+ D) \cdot C > 0$ for all curves $C$.

It follows that $K_{\mn}+ D$ is nef, and nef and big by induction.
Thus by the base point free theorem
(applied to the big and nef klt divisor $K_{\mn}+ D - \epsilon \Gamma$)
$m(K_{\mn}+ D)$ is basepoint free for $m \gg 0$. Since it
intersects every curve positively, it is thus ample. \qed\enddemo

The results above have some interesting geometric consequences:
\remark{3.7 Remarks} 
\roster
\item By (3.1.1) $\mm _{0,5}$ is isomorphic to $\pr 2.$ blown up at four 
points. Thus it is a del Pezzo surface of degree five. It is interesting to 
note that $K_{\mm _{0,5}}+D=-K_{\mm _{0,5}}$ is very ample and defines the 
anticanonical embedding of $\mm _{0,5}$ inside $\pr 5.$. For any $n$,
if $C$ 
is a vital curve, then $(K_{\mn}+D)\cdot C=1$. Thus it would be nice
to know if $K_{\mn}+D$ is very 
ample, for if it were, then under the corresponding map, 
every vital curve would be embedded as a line.
\item Note that $\lm _{0,5}$ is a log del Pezzo of rank one. It is easy 
to compute, using (3.2) that $\lm _{0,5}$ has two quotient singularities, 
one of index two and the other of index five. It is then easy, from the 
classification of log del Pezzos,
(see \cite{KeelMcKernan95}) to conclude that $\lm _{0,5}$ has one 
$A_1$-singularity and one singularity of type $(2,3)$. 
\item Note that the map $\map D_{T,T^c}.\ln.$
factors through $M_{|T| + 1,|T^c| + 1}/ S_{|T|} \times S_{|T^c|}$, but
not through the quotient by $S_{|T|+1} \times S_{|T^c| +1}$. Thus an
inductive study of $\overline{NE}_1(\ln)$ is problematic. In particular
one cannot obtain the analog of (3.9) as below.
\item By (3.1), $q_i^*(\ring .(1))$ is numerically equivalent to an effective 
divisor with support exactly $D$. It follows by (3.6) that for any curve 
$C \subset D$ there is some vital divisor which is negative on $C$. 
\endroster
\endremark 

\proclaim{3.8 Lemma} For any projective variety $T$, 
$N_1(\mn \times T) = N_1(\mn) \times N_1(T)$ under the map induced
by the two projections. The same map induces an isomorphism
$$
NE_1(\mn \times T) = NE_1(\mn) \times NE_1(T)
$$
\endproclaim
\demo{Proof} This follows from Theorem 2 of \cite{Keel92}.\qed\enddemo

\proclaim{3.9 Corollary} Fulton's conjecture for $NE^1(\mn)$ implies
the conjecture for $NE_1(\mn)$.
\endproclaim
\demo{Proof} Immediate from (2.2) and (3.8).\qed\enddemo

\demo{Proof of (1.2)} As we are going to use induction it is 
actually more convenient to prove a slightly stronger result. Let $M$ 
be any product of $\mm _{0,i}$. We will prove (1.2) for $M$. 
By a vital cycle on $M$ we mean a product of vital cycles on 
each component. We will continue to use the same notation, so for example by 
$D_{T}$ we mean the inverse image of this divisor from a projection onto one 
of the components of $M$.

Let $m$ be the dimension of $M$. When $m \leq 2$ it is easy to
check that vital curves generate the cone (see (3.7.1)).

$D$ has ample support by (3.6), and each
$D_{T,T^c}$ has anti-nef normal bundle by (4.5). 

When $m \leq 4$ then
we can apply (3.5),  and (2.6) inductively, to prove (1.2.3).

Let $R$ be an extremal
ray satisfying either (1.2.1) or (1.2.2). We show $R$ is spanned
by a vital curve by induction on $m$ which we may assume is at least
$5$. 

Note $M$ retracts onto any vital cycle, thus 
if $Z$ is any vital cycle, the restriction
$N^1(M) \rightarrow N^1(Z)$ is surjective. In particular
if $R$ is in the image of (the injection)
$i: \overline{NE_1}(Z) \rightarrow \overline{NE_1}(M)$,
then $R$ spans an extremal ray on $Z$.

Assume $R$ satisfies (1.2.2). By (2.6), $R$ spans an extremal
ray on some $D_T$. Since $D_T$ has anti-nef normal bundle, 
we can increase its
coefficient in $G$ to one, restrict to $D_T$, 
and apply adjunction and induction.

Assume $R$ satisfies (1.2.1). By (2.5), $R$ is spanned by a curve
$C \subset D_{T}$. By (3.7.4) we may assume $C \cdot D_T < 0$, so
$C$ does not deform away from $D_T$. 
By (2.7), $C \subset D_T$ is contracted by a map of relative
Picard number one, and so we can apply induction.
\qed \enddemo

\heading \S 4 Intersecting vital curves and divisors. \endheading

By a {\it marked point} of an $n$-pointed curve, we either mean one of 
the singular points of the curve, or one of the labeled points $\list p.n.$. 

{\it 4.1 Notation: Let $C$ be a vital curve. 
Let $G= G(C)$ be the $n$-pointed stable curve corresponding
to the generic point of $C$. $G$ has 
$n-3$ components, all but one of which contain
3 marked points, and exactly one of which contains $4$ marked
points. We call this last component $Q=Q(C)$, the distinguished component
of $G$. Let $s(C)$ be the number of singular points on $Q$, $l(C)$ be the 
number of labeled points. $C$ determines a decomposition of $\{1,2,\dots,n\}$
into $4$ disjoint subsets: $G \setminus Q$ has exactly $s(C)$ connected 
components. We decompose $\{1,2,\dots,n\}$ into those labeled points on
each of the components. Additionally we take the singleton
sets for each of the $l(C)$ labeled points on $G$. We call
this decomposition $P_C$.}

{\it There are $n-4$ singular points on $G$ (intersection points
of two components). Each singular $p \in G$ defines a decomposition,
by letting $T_p$ and $T^c_p$ be the labels on the two connected
components of $G \setminus \{p\}$. $C$ is the complete intersection
$\underset {p \in \sg(G)} \to \cap D_{T_p,T^c_p}$. Let
$A_{T_p}$ and $A_{T_p^c}$ be the connected components of 
$G \setminus \{p\}$. }

\proclaim{4.2 Lemma} Let $C$ be a vital curve.
\roster
\item For $p \in \sg(G)$
$$
\pi_{T_p}:\map D_{T_p,T_p^c}.M_{T_p \cup \{b\}}.
$$
contracts the vital curve $C$
iff $A_{T_p^c}$ contains the generic point of $Q(C)$. 
\item $q_i$ contracts $C$ iff $i$ is not one of the
labeled points of $Q(C)$ (in particular any $C$ with $l(C) =0$
is contracted).
\endroster
\endproclaim
\demo{Proof} (1) is immediate and (2) follows from (1) and 
(3.1.3). \qed\enddemo

\proclaim{4.3 Lemma} $P_C$ uniquely determines the
numerical class of $C$.  $K_{\mn}\cdot C = 2 - l(C)$. For any vital divisor
$D_{T,T^c}$ we have: 
\roster
\item $D_{T,T^c} \cdot C = -1$ iff $T$ or $T^c$ 
is one of the equivalence classes of $P_C$. Equivalently,
iff $T$ or $T^c$ is $T_p$ for some singular point $p \in Q$.
\item $D_{T,T^c} \cdot C =1$ iff $T$ or $T^c$ is the union of 
two equivalence classes.
\item Otherwise $D_{T,T^c} \cdot C =0$.
\endroster
\endproclaim
\demo{Proof} 
Since the vital divisors generate $\pic \mn.$ the description of
$D_{T,T^c} \cdot C$ implies the first statement. 
The expression for $K_{\mn}\cdot C$
follows from the expression for $D_{T,T^c} \cdot C$ using the
adjunction formula, since $C$ is a complete intersection of
vital divisors. 

Fix $p \in \sg(G)$ and let $S$ be the intersection
of the $D_{T_q,T^c_q}$ for $q \neq p$. Then $C \cdot D_{T_p,T^c_p}$
is the self intersection of $C$ in $S$. 
$S$ is a vital surface, and so it is either $\mm_{0,5}$ 
(which is $\pr 2.$ blown up in $4$ points)  or $\mm_{0,4} \times \mm_{0,4}$ 
(which is $\pr 1. \times \pr 1.$), and $C$ is a vital curve in $S$. In the
first case $C$ is a $-1$-curve, and in the second a fibre
of one of the two projections. Let $\gamma$ be the pointed stable
curve corresponding to a generic point of $S$. In the first case 
$\gamma$ has one component with $5$ marked points, and in the second
case, two components each with $4$ marked points. $G$ is obtained
as the limit as two of the marked points (on the same component) come
together at $p$. It's clear that the
first case occurs iff $p \in Q$, whence (1). Note the
argument shows that if $C \subset D_{T,T^c}$ then $C \cdot D_{T,T^c}$
is either $0$ or $1$.

If $D_{T,T^c} \cdot C > 0$ then $D_{T,T^c} \cap C$
is a vital point of $C =\mm_{0,4}$,  i.e. a reduced point, thus
$D_{T,T^c} \cdot C =1$. This occurs if $T$ or $T^c$ is a union
of two equivalence classes of $P_C$ and every 
vital divisor of $C$ can be obtained in this way. Since each vital
cycle is uniquely a complete intersection of vital divisors, this gives
(2). 

Since the possibilities with $C \cdot D_{T,T^c}$ nonzero are 
classified by (1) and (2),  (3) follows. \qed\enddemo

\proclaim{4.4 Corollary} The numerical class of $\t{C}$ is determined
by the cardinalities of the subsets in $P_C$. If these cardinalities
are $a$, $b$, $c$, $d$ then
$$
C \cdot \sum r_i B_i = -r_a  - r_b - r_c - r_d + r_{a + b} + 
r_{a + c} + r_{a + d}
$$
where we define $r_1 = 0$ and $r_i = r_{n-i}$
for $i > [n/2]$. \endproclaim

\proclaim{4.5 Lemma} 
$$ 
N_{D_{T,T^c}}\mn = (q_b \circ \pi_{T})^*(\ring .(-1)) \otimes
(q_b \circ \pi_{T^c})^*(\ring .(-1)) .
$$
\endproclaim
\demo{Proof} Since the vital curves generate $N_1$ 
we only need to check how both 
sides intersect a vital curve $C \subset D_{T,T^c}$. By (3.1.2),
and (4.3) the possible values of these intersections are $0$ and
$-1$, and it is enough to show 
show $D_{T,T^c} \cdot C = -1$ iff one of the two maps
$q_b \circ \pi_T$ or $q_b \circ \pi_{T^c}$ fails to contract
$C$. 

By (4.2.1) we may assume that $\pi_T$ is finite on $C$ (otherwise
switch $T$ and $T^c$).
By (4.2.1) and (4.3.1), $D_{T,T^c} \cdot C = -1$ iff $b$ is a labeled point
of $Q(\pi_T(C))$, thus by (4.2.2), iff $q_b$ is finite on
$\pi_T(C)$. \qed \enddemo

\proclaim{4.6 Lemma} Let $C$ be a vital curve, and let 
$$
D_C = \underset {p \in \sg(G) \cap Q} \to \sum D_{T_p,T^C_p}. 
$$
$(K_{\mn}+ D_C) \cdot C = -2$. $K_{\mn}+ D + 1/s D_C$ intersects vital
curves non-negatively, and vanishes on exactly those vital curves 
numerically equivalent to $C$. 
\endproclaim
\demo{Proof} Immediate from (4.3). \qed \enddemo

\remark{4.6.1 Remark} (1.1) and the basepoint free theorem imply
$K + D + 1/s D_C$ is eventually free, and thus $C$ spans an extremal
ray. Presumably this could be checked directly. \endremark

 The following is immediate:

\proclaim{4.7 Lemma} Let $T \subset \{1,2,\dots,n\}$ with
$|T| \geq 3,|T^c| \geq 2$. For $i \in T$, 
$$
D_{T \setminus \{i\},T^c \cup \{i\}}|_{D_{T,T^c}} =
D_{ib,T \setminus \{i\}} \times M_{T^c \cup \{b\}}
$$
under the canonical product decomposition. There is no other vital divisor 
with the same restriction.
\endproclaim

\proclaim{4.8 Lemma} Suppose there is a numerical equality
$$
\sum_{i=2}^{k}  m_i B_i \thicksim F
$$
and either $F$ is nef, or both sides are effective and have no
divisor common to their supports. Then
$$
\aligned
rm_{r-1} &\geq (r-2) m_r \\
(n-r) m_{r+1} &\geq (n-r-2) m_r 
\endaligned
\qquad
\aligned 
\text{ for } &3 \leq r \leq k \\
\text { for } &2 \leq r \leq k -1.
\endaligned
$$
When the left hand side is effective, it is either
trivial, or has support exactly $D$. (1.3.1-4) hold. 
\endproclaim
\demo{Proof} We prove the first inequality, the argument for the second is 
analogous. The final remarks follow from the inequalities.

Choose $T$ with $|T| =r$. Let $Z_r$ be the general fibre of 
$$
\map M_{T \cup \{b\}}.M_{T}..
$$
Let $p \in M_{T^c \cup \{b\}}$ be a general point, and let
$D_{T,T^c} \supset C_r = Z_r \times \{p\}$. By (3.1), (4.5) and
(4.7) we have
$$
C_r \cdot B_i = \left\{ \aligned r &\text{ if } i = r-1 \\
                                 -(r -2) &\text{ if } i =r \\
                                 0 &\text{ otherwise } \endaligned \right.
$$
The inequality is obtained by intersecting both sides with $C_r$.

 Note that (1.3.1-3) follow immediately and that (1.3.4) then follows from 
(2.3).\qed\enddemo

\heading \S 5 $NE_1(\ln)$ for small $n$  \endheading

Given (4.8), it is natural to hope that every nef divisor on $\ln$ is 
eventually free. The obvious approach is to try to use the basepoint free 
theorem, and thus to realise some positive multiple of a big
nef class $E$ (pulled back from $\ln$) as a klt divisor $K_{\mn}+ \Delta$. 

\proclaim{5.1 Lemma} If $E$ is a big nef class on a normal
$\Bbb Q$-factorial variety $M$, and there is a divisor $\Delta$ with
$K_M+\Delta$ klt and numerically equivalent to a positive 
multiple of $E$, then the extremal subcone of $\overline{NE}_1(M)$
supported by $E$ is rational polyhedral, and is contracted by a
log Mori fibre space. If $M = \mn $, the
subcone supported by $E$ is spanned by vital curves. \endproclaim
\demo{Proof} By (2.2) we have $E = A + Z$ where $A$ is ample
and $Z$ is effective.  If $V \subset \overline{NE}_1(M)$ is
the extremal subcone supported by $E$, then 
$K_M+ \Delta + \epsilon Z$ is negative on $V \setminus 0$. 
Thus the result follows from the cone and 
contraction theorems, together with (1.2) 
\qed \enddemo

Let $E$ be a nef divisor on $\mn$, pulled back from $\ln$.

By (3.5), for $n \leq 7$ the conditions of (5.1) are satisfied
(if $K + \Gamma$ is klt and trivial, let $\Delta = \Gamma + \epsilon E$).

In general, 
by (3.5) and (4.8), replacing $E$ by a large multiple
one has $E = K_{\mn}+ \Delta$ for some $\Delta$ supported on
$D$. We can try to make $\Delta$ a boundary by subtracting off
part of $E$, thus we are lead to consider:

\proclaim{5.2 Definition-Lemma} Let $E$ be a non-trivial nef
class on $\mn$, pulled back from $\ln$ with $n \geq 8$.
Then there is a unique effective class $\Delta_E$ with
the following properties
\roster
\item $\Delta_E$ has support a proper subset of $D$
\item $K_{\mn}+ \Delta_E = \lambda E$ for some $\lambda > 0$.
\endroster
\endproclaim
\demo{Proof} For any $\lambda$, $-K_{\mn}+ \lambda E$ is pulled back
from $\ln$, thus by (4.8), (1) is the requirement that $\Delta_E$ be on the
boundary of $NE^1$. Since $E$ is in the interior of $NE^1$, and
by (3.5), $-K_{\mn } \not \in NE^1$, the result is clear.
\qed \enddemo 

Notation: {\it For the next corollary, define the integer function
$f(a,b,c,d)$ to be $2$ minus the number of variables equal
to one.

 We will say that $P_n$ holds if for a given integer
$n$ the following implication holds:

Let $\list r.n-1.$ be a collection of non-negative real numbers,
with $r_1 = 0$, $r_i = r_{n-i}$, and $r_j =0$ for some $2 \leq j \leq k$.
If
$$
f(a,b,c,d) + r_{a+b} + r_{a + c} + r_{a + d} \geq r_a + r_b + r_c + r_d
$$
for every set of positive integers $a$, $b$, $c$, $d$ with
$n = a + b + c +d$, then $r_i < 1$ for all $i$.}

\proclaim{5.3 Corollary} $\Delta_E$ is a pure boundary for
every non-trivial nef class pulled back from $\ln $ iff
$P_n$ holds. \endproclaim
\demo{Proof} By (4.8) $P_n$ is equivalent to the statement:
If $\sum r_i B_i$ has support a proper subset of $D$ and 
$(K_{\mn}+ \sum r_i B_i) \cdot C \geq 0$ for all vital curves $C$, then 
$\sum r_i B_i$ is a pure boundary. Thus the only thing to
show is that if $\Delta_E$ is a pure boundary for every
non-trivial nef class, then the images of vital curves generate
$NE_1(\ln)$. This follows from (5.1). \qed \enddemo

For a given $n$ it is straightforward to check whether or not $P_n$ holds:

\proclaim{5.4 Lemma} $P_n$ holds for $8 \leq n \leq 11$.
\endproclaim
\demo{Proof} We will check $P_9$. The cases $n=8$, $10$ and $11$ are similarly
checked. 

Let $\list r.8.$ be a collection of non negative numbers,
as in the definition of $P_9$. From the sums 
$$
\align 
1+2+3+3 &=9\\
1+1+1+6 &=9\\
1+2+2+5 &=9
\endalign
$$
we obtain the inequalities 
$$
\align
1 + r_4 & \geq r_3 + r_2 \\
2r_3 & \geq r_4 \\
3r_2 & \geq r_3 + 1
\endalign
$$

The result follows easily by considering in turn the possibilities $r_4=0$, 
$r_3=0$ and $r_2 =0$. \qed \enddemo

Observe that (1.3) follows from (5.1), (5.3) and (5.4).

We have checked that $P_n$ holds for several $n \geq 12$, and we
suspect that (if motivated) one could prove this for all such $n$.

\end